\documentclass[12pt]{amsart}

\usepackage{amsmath}
\usepackage{amssymb}
\usepackage{graphicx}
\usepackage{epsfig}
\usepackage{geometry}
\usepackage{amsfonts}
\usepackage{graphics}

\newtheorem{theorem}{Theorem}[section]

\newtheorem{proposition}[theorem]{Proposition}

\begin{document}
\parskip0.5em

\title {Genus of Embedded Graphs in Orientable Closed Surfaces}

\author {Lorena Armas-Sanabria}
\address{CONAHCYT-Universidad Aut\'onoma del Estado de Morelos, Av. Universidad 1001, Cuernavaca 62209, Mor., M\'exico.}
\email{lorenaarmas089@gmail.com}

\author{V\'ictor N\'u\~nez}
\address{Centro de Investigaci\'on en Matem\'aticas, A.P.~420, Guanajuato 36000, Gto., M\'exico.}
\email{victor@cimat.mx}

\thanks{}



\begin{abstract}

 We give an algorithm to calculate the minimal and maximal genus of
 the orientable closed surface where a graph $G$ can be embedded. For
 this, we construct some special branched coverings of the
 2-sphere. We apply this algorithm to calculate the
 orientable genus and maximal genus of  some Snarks graphs.  
\end{abstract}

\maketitle

\section{Introduction}

The \emph{orientable genus} of a graph $G$ is defined as the minimal genus of
an orientable closed  surface where the graph can be embedded. The
\emph{maximal orientable genus} of $G$ is the maximal 
genus  of an orientable closed surface where the graph can be
cellularly  embedded, and remember that an embedding is cellular if the complementary regions of the graph in the
surface consist of open disks. The set of numbers bounded by the orientable genus and the maximal orientable genus is called
the genus range of $G$.

Given a graph $G(V,E)$ where $V$ is the set of vertices and $E$ is the
set  of edges of $G$, there is a way to calculate the
genus of $G$, by using what is called a rotation system \cite{R}. A
rotation system consists in an assignment of a numeration  to the edges and
vertices, in such a way that each vertex has all its incident edges
numerated,  and then a cyclic permutation is assigned  by considering
the incident edges to each vertex $v_i$. It is not difficult to see that this
cyclic permutation defines an embedding of $G$ on an orientable closed
surface. This fact was first observed by Edmonds \cite{E}. There are many possibilities of assigning a cyclic
permutation to each vertex,  these possibilities define all the
possible embeddings of the graph in different surfaces. 
The first algorithm to calculate the genus of a graph was given by Youngs \cite{Y},
which consists in producing all rotation systems, and for each such system
determine the cycles who are to bound a disk in the surface. Then knowing
the number of vertices, edges and cycles that bound disks, the genus of the
embedding can be calculated by means of the Euler characteristic.
This algorithm determines all the genus range of $G$, but it requires many steps, for if $deg(v_i)$ denotes the degree of vertex $v_i$,
note that there are $\prod _{v\in V(G)}(deg(v_i) - 1)!$ rotation systemss. Also, by a result of Duke \cite{D}, it is known that if $k$ 
is a number in hte renge genus of $G$, then there exists a celular embedding of $G$ in a closed orientable surface of genus $k$.

There are many other methods to determine the genus of a graph,
see the survey paper of Stahl \cite{S}.

The problem of calculating the genus of a graph is an
$NP$-complete problem \cite{T}. 
There are some bounds for these genera \cite{R}, which usually, are not very
sharp. A way to show that a graph cannot be embedded in some surface,
is determining prohibited minors of graphs for a given orientable, 
closed surface $F$;  but it is  also  a difficult problem to find 
complete sets of prohibited minors. Given
an arbitrary  orientable closed surface $F$, a
complete  set of prohibited minors is not known, except for the case  $F$ is
a 2-sphere, which is the Kuratowski Theorem, which states that a given graph
is planar if and only if it does not have a subgraph homeomorphic to
$K_5$,  or to $K_{3,3}$. Also there are some known
prohibited minors for the case of the torus, and in the non-orientable
case, a complete set of prohibited minors  for the
projective plane is  known \cite{A}. 

 There are several algorithms for embedding  graphs of arbitrary genus  given by Filotti, Miller and Reif \cite{FMR}, and by 
Djidjev and Reif \cite{DR}, but in \cite{MK}, it is pointed out that they are incorrect and that there is not way  to fix them without creating algorithms which 
take exponential time.
Mohar, in \cite{M}, \cite{M1},  proved that if $S$ is a fixed surface then there is a linear time algorithm  to 
know if an arbitrary  graph  $G$ embeds in $S$ (and finds an embedding) or, finds a subgraph of $G$ which is a subdivision of some graph in the set of forbidden graphs of $S$ (see also \cite{MT}). However, in \cite{MK}, it is claimed that this algorithm is very complex and very difficult to implement.

  In \cite{G}, it is constructed a quadratic-time algorithm for calculating the genus distribution of the graphs in a family ${\mathbb{F}}$, consisting of graphs of fixed treewidth and bounded degree, which complements an algorithm of  Kawarabayashi, Mohar and Reed \cite{KMR}  for
  calculating the minimum genus of a graph of bounded treewidth in linear time.

In this paper we develop an algorithm to calculate the genus and
the maximal genus (in fact, the genus range)
of a graph, which is equivalent to the construction of the rotation
systems of the graph, but in a nicer setting. 
Our main technique is the construction of certain branched coverings
of the 2-sphere. This gives short and clear proofs of the correctness
of the algorithm, and also explicit
embeddings.

The paper is organized as follows. In Section~\ref{sec2} we recall the
elements of branched covering theory.  
In Section~\ref{sec3} we
describe the algorithm to calculate the orientable genera of a
graph $G$, and we show that the algorithm works. In Section~\ref{sec4} we
work an example, and finally in Section~\ref{sec5} there is a table with
the orientable genus and maximal orientable genus for some
Snarks.

\section{Branched covers of the 2-sphere} 
\label{sec2}
For a nice account of branched covering theory, see \cite{BE}.
 A finite-to-one, open and proper map $\varphi:M\rightarrow N$ between
 two $n$-manifolds  is called
a \emph{branched covering}.
The usual way to check that an open map $\varphi:M\rightarrow 
N$ is a branched covering, is to find
a codimension 2 properly embedded submanifold $k\subset N$  such that
$\varphi:M-\varphi^{-1}(k)\rightarrow N-k$
is a $d$-fold covering space.
One says that \emph{$\varphi$ is branched along $k$}, and that $k$ is the
\emph{branching} of $\varphi$.
We write $\varphi:M\rightarrow (N,
k)$ for a branched covering $\varphi:M\rightarrow N$ with branching~$k$.

For a given $d$-fold  branched covering $\varphi:M\rightarrow (N, k)$,
there is an associated
\emph{representation} (that is, a homomorphism) of the fundamental group
$\pi_1(N-k)$ into~$\Sigma_d$, the symmetric group on $d$ symbols,
$\omega_\varphi:\pi_1(N-k)\rightarrow \Sigma_d$, as follows: Fix a base point~$x_0\in N-k$, and
number $\varphi^{-1}(x_0)=\{x_1,\dots,x_d\}$. For $[\alpha]\in
\pi_1(N-k,x_0)$, consider the lifting of $\alpha$ at the point $x_i$
which ends, say, in the point $x_j$. Then define~$\omega_\varphi([\alpha])(i)=j$.

For a given representation $\omega:\pi_1(N-k)\rightarrow \Sigma_d$,
one has an associated $d$-fold branched covering
$\varphi_\omega:M\rightarrow (N, k)$: it is the completion of the
covering space~$M_0\rightarrow N-k$ corresponding to the subgroup
$\omega^{-1}(\Sigma_{d-1})$. 

Equivalence classes of $d$-fold connected branched covers
$M\rightarrow (N,k)$  are in
one-to-one correspondence with conjugacy classes of transitive
homomorphisms $\pi_1(N-k) \rightarrow \Sigma_d$.

\begin{figure}

\includegraphics[angle=0, width=7true cm]{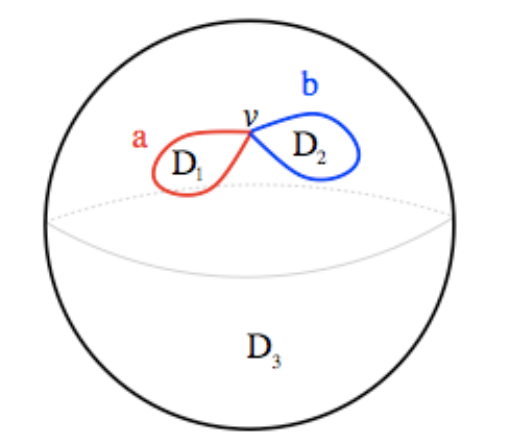}
\caption{A cellular decomposition of the 2-sphere.}
\label{fig1}

\end{figure}

We are interested in branched coverings of the 2-sphere $S^2$ branched
along three points $x_1,x_2,x_3\in S^2$. These coverings are related
to a  cell decomposition of the~2-sphere consisting of one vertex $v$,
two edges $a$ and $b$, and three disks
$D_1,D_2,D_3$ oriented in such a way that  $\partial
D_1=a, \partial D_2=b, \partial D_3 =(a*b)^{-1}$. We regard $x_i$ as
the center of~$D_i$ (see Figure~\ref{fig1}).

The fundamental group of the punctured 2-sphere is free of rank two and
admits a presentation  $\pi_1(S^2-\{x_1,x_2,x_3\})\cong\langle
\bar x_1,\bar x_2,\bar x_3:\bar x_1\bar x_2\bar x_3=1\rangle$, where the symbol $\bar x_i$ is the class
of a closed curve on $S^2$ starting at the vertex $v$ and encircling
the point~$x_i$. Alternatively, using the cell decomposition of $S^2$
above, we 
obtain a presentation $\pi_1(S^2-\{x_1,x_2,x_3\})\cong\langle
a,b,c:abc=1\rangle$ where $a\simeq\bar x_1$, $b\simeq\bar x_2$, and
$c=(a*b)^{-1}\simeq\bar x_3$. Therefore an $d$-fold branched covering of $S^2$
branched along $\{x_1,x_2,x_3\}$ is determined by a triple of
permutations $\alpha,\beta,\gamma\in \Sigma_d$, such that
$\alpha\beta\gamma=1$. The triple $(\alpha,\beta,\gamma)$
is called a \emph{constellation} in \cite{LZ}. An arbitrary pair of permutations $\alpha,\beta\in
\Sigma_d$ gives a constellation, $(\alpha,\beta,\gamma)$, setting $\gamma=(\alpha\beta)^{-1}$.

A constellation $\alpha,\beta,\gamma\in \Sigma_d$ defines a representation
$\omega: \pi _1 ( S^2 - \{x_1,x_2,x_3\}) \rightarrow \Sigma_d$ such that
$\omega(a)=\alpha$, $\omega(b)=\beta$, and $\omega(c)=\gamma$. The
associated branched covering~$\varphi_\omega:F\rightarrow S^2$ can be
constructed in two steps. 
 First, construct the covering space of
$a\vee b$ induced by  the restriction $\omega:\pi_1( a\vee
b)\rightarrow \Sigma_d$, and, secondly,
construct the branched coverings of the disks $D_i$ branched along
$x_i$ associated to the restrictions~$\omega: \pi _1 ( D_i - \{x_i\})
\rightarrow 
\Sigma_d$. These covers have as many components as the number of cycles in
a disjoint cycle decomposition of
$\omega(\partial(D_i))$. Finally assemble together these coverings
into $F$ with liftings of  
the glueing maps~$\partial D_i\rightarrow a\vee b$.
  
The $d$-fold branched cover
$\phi : F\rightarrow S^2$ gets a cell decomposition for $F$ induced by the
cell decomposition of $S^2$,  with  vertices $\varphi^{-1}(v)$, edges
$\varphi^{-1}(a\cup b)$, and 2-cells~$\varphi^{-1}(D_1\cup D_2\cup D_3)$.

\section {Orientable genera of graphs} 
\label{sec3}
In this section we describe an algorithm  to calculate the orientable
genera of a graph $G$. For a set $Y$, we write $\Sigma(Y)$ for the symmetric
group in the symbols in~$Y$. If $\#(Y)=k$, we write $\Sigma_k=\Sigma(Y)$. If
$\gamma\in \Sigma_k$, $|\gamma|$ is the number of cycles in a
decomposition of $\gamma$ into disjoint cycles in $\Sigma_k$ (including
trivial cycles for the fixed points of $\gamma$).
 
\subsection {The algorithm}

\vskip 15pt

 We start with a finite graph $G$ with vertices
 $v_1,\cdots, v_n$, and edges
$f_1,\cdots, f_e$. Without loss of generality, we assume $G$ simple
and with all vertices of degree at least 3.

\begin{figure}

{\includegraphics[angle=0, width=6true cm]{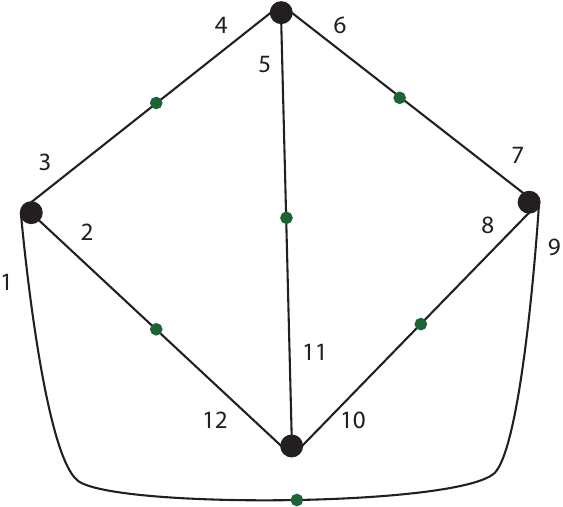}}
\caption{$K_4$ with darts.}
\label{fig2}
\end{figure}

\noindent {\bf 1 )} Add new fake vertices
$w_1,\cdots, w_e$ where each vertex  $w_i$ is the middle point of
$f_i$. We obtain two new edges for each $f_i$, which are called
\emph{darts} in \cite{LZ} (See Figure \ref{fig2}).

\noindent {\bf 2 )} Label each dart with a number in $\{ 1,2,\cdots ,2e \}$, and define the
permutation $\alpha = (d_1,{d'}_1)(d_2,{d'}_2) \cdots (d_e,{d'}_e)\in S_{2e}$, where $d_i, {d'}_i$
are the numbers of the darts of the original edge $f_i$. 

\noindent {\bf 3 )} For each vertex $v_i$, consider $c_i$, the conjugacy class of the $t_i$-cycle
$(s_{i,1},s_{i,2},\dots , s_{i,t_i})$ in ${\Sigma(\{s_{i,1},s_{i,2},\dots , s_{i,t_i}\})}$, where
$s_{i,1},s_{i,2},\dots , s_{i,t_i}$ are the numbers of the darts
incident in $v_i$. 

\noindent {\bf 4 )} For each choice $\beta_1\in c_1,\beta_2\in c_2,\dots,\beta_n\in
c_n$, define the permutation
$\beta = \beta_1\beta_2\cdots\beta_n$.

\noindent {\bf 5 )}  Collect the numbers
$(2-(|\beta|-|\alpha|+|(\alpha\beta)^{-1}|))/2$ in a set $X$, for all possible $\alpha, \beta$.

Note that each element of a set $c_i$ is a permutation system around the vertex $v_i$, and a choice of permutation
$\beta$ is just a choice of a rotation system for the graph.

Let $g(G)$ and $g_M (G)$ be the orientable genus and the maximal orientable genus of $G$, respectively. 

\vskip 15pt

\begin{proposition} Let $G(V,E)$ be a  graph such that $deg(v_i)\geq3$ for all $i$. Let  $X = \{ (2-(|\beta|-|\alpha|+|(\alpha\beta)^{-1}|))/2 \hskip 5pt  | \hskip 5pt$ for all possible  $ \alpha, \beta \}$. Then $g(G)= \min(X)$ and $g_M(G) = \max(X)$.
\end{proposition}

\subsection{A branched covering}

We verify that the described algorithm works. For that we show that
the genera of all orientable surfaces $F$ such that $G$ embeds in $F$,
are collected in the set $X$ above.

\begin{proposition}Let $G(V,E)$ be a graph such that $deg(v_i) \geq3$ for all $i$. Then the genus range of $G$ is given in the set $X$.
\end{proposition}

\noindent {\it Proof}.
We consider the graph $G$ with its edges subdivided in darts as before.
Assume the graph $G$ is cellularly embedded in an oriented surface $F$, and
let~$U$ be a regular neighborhood of $G$ in $F$. 

Inside $U$, substitute each vertex $v_i$ with a fat vertex which is a polygon of
$t_i$ vertices
where  $t_i$ is the valency of the vertex  $v_i$. We
visualize this fat vertex as centered in~$v_i$, and intersecting the
graph $G$ in $t_i$ points of the darts incident in $v_i$. As an example of this construction, we can see Figure \ref{fig3}.

The vertices of this  $t_i$-gon are labeled with the numbers of the
incident darts  in~$v_i$. 

Also, inside $U$, substitute each edge $f_i$ 
of $G$ with two parallel edges whose ends are the intersection
points of $f_i$ with the polygons corresponding to the ends of $f_i$.

We obtain a graph $G'$ inside $U$, whose vertices are the vertices of
all of the
$t_i$-gons, and whose edges are the edges of the  
$t_i$-gons and the added parallel edges, one pair for each original
edge of $G$. 

This graph $G'$ defines a pair of permutations,
$\alpha = (d_1,d_1')(d_2,d_2') \cdots (d_e,d_e')$
where $d_i, d_i'$ are the ends of the parallel edges corresponding to
the original edge $f_i$,
and
$\beta=(s_{1,1},s_{1,2},\dots,s_{1,t_1})(s_{2,1},s_{2,2},\dots,s_{2,t_2})\cdots
(s_{n,1},s_{n,2},\dots,s_{n,t_n})$ where $s_{i,1},s_{i,2},\dots ,
s_{i,t_i}$ are the numbers of the vertices of the 
$t_i$-gon corresponding to the original vertex $v_i$ and are cyclically
ordered by the local orientation of $F$ at $v_i$.

This determines a $(2e)$-fold covering space
$\psi:\tilde{G} \rightarrow a\vee b$ corresponding to the
representation $\pi_1(a\vee b)\rightarrow \Sigma_{2e}$ such that
$a\mapsto\alpha$, and $b\mapsto\beta$.
Notice that in this covering,
the preimage of  $a$
has two edges for each transposition $(d_i,d_i')$ of $\alpha$, and the
preimage of $b$ has a $t_i$-gon for each cycle $(s_{i,1},s_{i,2},\cdots , s_{i,t_i})$ of $\beta$.
That is,  $\tilde{G} = G'$.

We can extend $\psi$ into a $(2e)$-fold branched covering
$\varphi:F\rightarrow(S^2,\{x_1,x_2,x_3\})$ corresponding to the
representation
$\pi_1(S^2-\{x_1,x_2,x_3\})\cong\langle a,b,c:abc=1\rangle\rightarrow \Sigma_{2e}$ such that $a\mapsto\alpha$,
$b\mapsto\beta$, and $c\mapsto \gamma=(\alpha \beta)^{-1}$.

For a permutation $\eta\in\Sigma_n$, we write $o(\eta)$ for the order
of $\eta$ in the group $\Sigma_n$.
Write $\alpha=\alpha_1\cdots\alpha_e$, $\beta=\beta_1\cdots\beta_n$,
and $\gamma=\gamma_1\cdots\gamma_m$ for disjoint cycle decompositions
of $\alpha,\beta$ and $\gamma$ in $\Sigma_{2e}$. Then
$\alpha_i=(d_i,d'_i)$, and $\beta_i=(s_{i,1},s_{i,2},\cdots ,
s_{i,t_i})$. Let $\chi (F)$ be the Euler characteristic of $F$. Since $F$ is obtained as a branched covering of $S^2$, by the Riemann-Hurwitz formula we have:
\begin{displaymath}
\begin{array}{rl}
\displaystyle
\chi(F)&\displaystyle=2e\chi(S^2)-\sum_{i=1}^e(o(\alpha_i)-1)-\sum_{j=1}^n(o(\beta_j)-1)-\sum_{k=1}^m(o(\gamma_k)-1))\\
&\displaystyle=4e-((\sum_{i=1}^e2)-e)-((\sum_{j=1}^n t_{j})-n)-((\sum_{k=1}^m o(\gamma_k))-|\gamma|)\\
&=4e-(2e-e)-(2e-n)-(2e-|\gamma|)\\
&=n-e+|\gamma|\\
&=|\beta|-|\alpha|+|(\alpha\beta)^{-1}|.\\
\end{array}
\end{displaymath}

Since $\chi(F)=2-2g(F)$, we see that
$$g(F)=\frac{2-(|\beta|-|\alpha|+|(\alpha\beta)^{-1}|)}{2},$$
as needed.

Notice that, if we change the orientation of $F$, the permutations
$\alpha$, $\beta$ and $\gamma$ constructed, are just the inverses of
the permutations of the text, and, therefore,
have the same cycle
structure in $\Sigma_{2e}$, giving the same formula for $g(F)$.

Also notice that, if  $g$ is an arc in the 2-sphere connecting
$x_1$ and $x_2$ and intersecting $a\cup b$ exactly in the vertex $v$
(see Figures~\ref{fig1}, \ref{fig3}), then $\varphi^{-1}(g)=G$.

On the other hand, if we are given a pair of permutations $\alpha$ and $\beta$, we can construct a surface $F$ 
containing $G$, such that the corresponding associated permutations are precisely $\alpha$ and $\beta$.
To see this, consider an oriented polygon $t_i$ for each vertex $v_i$, which carries the permutation $\beta$ restricted to $v_i$. For an edge $f_k$ connecting vertices $v_i$ and $v_j$,  
attach a rectangle connecting the polygons $t_i$ and $t_j$, which preserves the orientation of the polygons. We can assume that the union of the disks and the band
contains the edge $f_k$. Doing this for all edges, produces a surface with boundary, and then by attaching disks to each of its
boundary components we get a closed surface $F$ in which $G$ embeds cellularly. 
\qed

Now, just see that Proposition 3.1 is a corollary of Proposition 3.2.

\section{Example} 
\label{sec4}

 In this Section, we 
apply  the algorithm to $K_4$, the complete graph of $4$ vertices as in Figure~\ref{fig2}.

 From Figure~\ref{fig2}, we obtain the permutation
$\alpha=(1,9)(2,12)(3,4)(5,11)(6,7)(8,10)$. We choose, as in the
figure, $\beta=(1,2,3)(4,5,6)(7,8,9)(10,11,12)$, and then $\gamma=(\alpha\beta)^{-1}= (1,4,7)(2,9,10)(3,12,5)(6,11,8)$.
Here $e=6$, $t_j=3$ for $j=1,2,3,4$, and $l_k=3$ for $k=1,2,3,4$.

Using the  Riemann-Hurwitz formula we obtain: 

$$g(\tilde X) = 2-(|\beta|-|\alpha|+|(\alpha\beta)^{-1}|)/2=2-(4-6+4)/2=0$$

So, $\tilde X$ is the $2$-sphere. But we did know that $K_4$ is planar.

\begin{figure}

{\includegraphics[angle=0, width=12true cm]{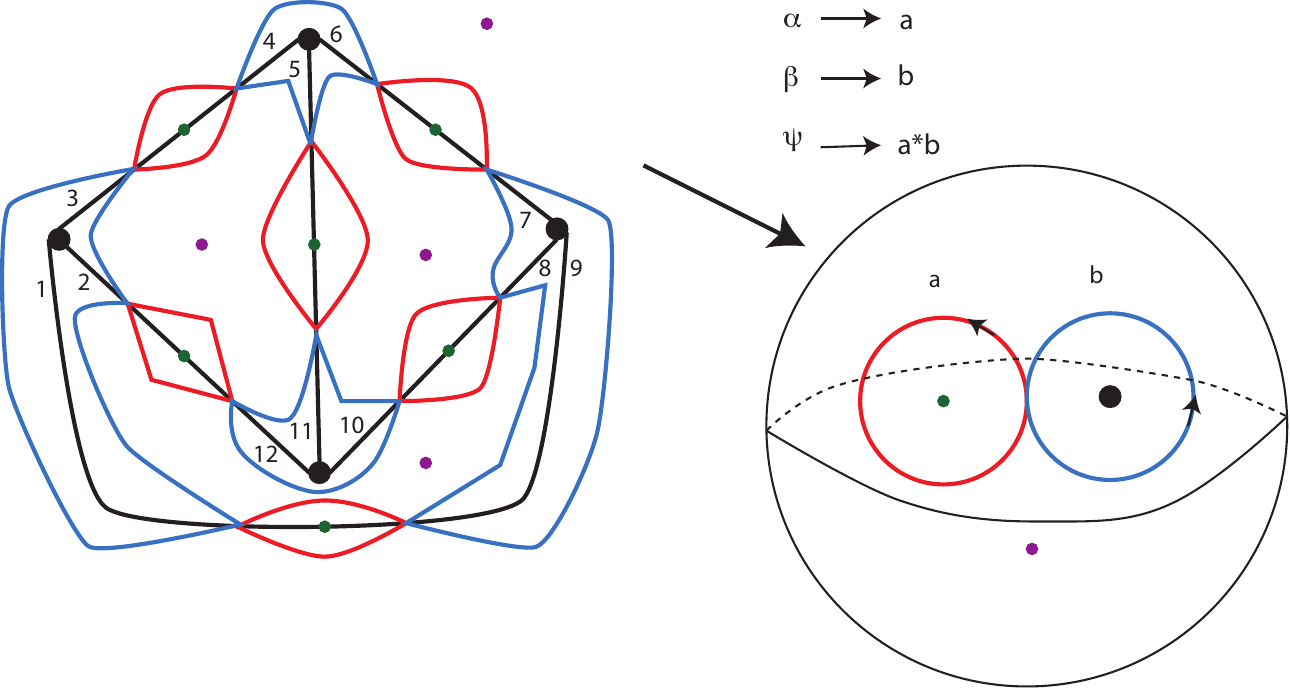}}
\caption{The 2-sphere as a branched cover of the 2-sphere.}
\label{fig3}
\end{figure}

The graph $G'$ looks like the graph in the left side of Figure \ref{fig3}, and in that figure one can see the corresponding branched cover of the 2-sphere.

Now $g_M(K_4) = 1$, and this genus can be achieved with the representation
$\alpha = (1,9)(2, 12)(3, 4)(5, 11)(6, 7)(8, 10) $ as before, but 
$\beta = (1, 2, 3)(4, 5, 6)(7, 8, 9)(10, 12, 11)$. Then $(\alpha \beta)^{-1} = (4,7,1)(10,5,3,12,8,6,11,2,9)$,
and using the genus formula again we obtain:

$$g(\tilde X) = (2-(|\beta|-|\alpha|+|(\alpha\beta)^{-1}|))/2=(2-(4-6+2))/2=1$$

\section { Orientable genus and maximal orientable genus of some Snarks}  
\label{sec5}
\medskip
 
In this section, we list the orientable genus and maximal orientable
genus of some Snarks. Recall that a Snark is a 3-regular graph, which is 4-edge connected
and it is not 3-edge-colourable. Many examples of Snarks are known, this is an interesting
but mysterious class of graphs. We use the list of snarks and the notation given in
\cite{RW}.  A  program for 3-regular graphs implementing the algorithm
of Section~\ref{sec3}, and written in the GAP
programming language is available from the authors on request.

\centerline {\textbf {}}

\begin{center}

 \begin{tabular}{|c|c|c|c|c|c|c|}

\hline

  $Graph$ & Orientable genus& Maximal or. genus   \\

\hline 

Petersen  & 1   &  3  \\
\hline

Tietze &1    &3  \\

\hline 

 Blanusa 1 & 2      &5  \\

\hline 

 Blanusa 2 & 1   &  5   \\

\hline 

 j5 &2    &5    \\

\hline 

Sn13  & 2   & 6   \\

\hline 

 Sn18 &3   & 6    \\

\hline 

 Sn19 & 2   & 6  \\

\hline 

Sn20 &2  &6   \\

\hline 

Sn21 &2 &6   \\

\hline 

 Sn22& 2   & 6    \\

\hline 

Sn23 & 2   & 6    \\

\hline 

Sn24 & 2  &6    \\

\hline 

Sn25 & 2   & 6   \\

\hline 

Sn26 &2   &6  \\

\hline 

Sn27 & 2   & 6  \\

\hline 

Sn28  & 2   &6   \\

\hline 

 Sn29 &  2  &6    \\

\hline 

 Sn27bis &2  & 6   \\

\hline 

Sn28bis & 2  & 6   \\

\hline 

 CelminsSwart1& 2   & 7    \\

\hline 

LoupequineTypeI  &  1  & 7  \\

\hline 

Type2&  2 &  7   \\

\hline 

CelminsSwart2&   2 &   7  \\

\hline
FlowerJ7 & 2 &7  \\

\hline
 doubleStar & 3 & 8 \\

\hline Double St  &  3  &  8 \\

\hline
 Type1 & 1 & 9\\

\hline Type2 34  &  2  &  9 \\

\hline
FlowerJ9 & 2 &9  \\

\hline
\end{tabular}

\end{center}

\textbf{Acknowledgement.} The first author is supported by a fellowship \textit{Investigadoras por M\'exico} from CONAHCYT. Research partially supported by the grant PAPIIT-UNAM IN117423.

\end{document}